\documentstyle[12pt,amscd,amssymb]{amsart}
\theoremstyle{plain}
 \newtheorem{thm}{Theorem}[section]
 \newtheorem{lem}[thm]{Lemma}
 \newtheorem{prop}[thm]{Proposition}
 \newtheorem{cor}[thm]{Corollary}
 
\theoremstyle{definition}
 \newtheorem{defn}{Definition}[section]
\theoremstyle{remark}
 \newtheorem{rem}{Remark}[section]

\newcommand{ \Supp}{\operatorname{Supp}}
\newcommand{\Ext}{\operatorname{Ext}}
\newcommand{\Hom}{\operatorname{Hom}}
\newcommand{\codim}{\operatorname{codim}}

\newcommand{\Aut}{\operatorname{Aut}}
\newcommand{\rk}{\operatorname{rk}}

\newcommand{\NS}{\operatorname{NS}}
\newcommand{\coker}{\operatorname{coker}}
\newcommand{\Pic}{\operatorname{Pic}}
\newcommand{\ch}{\operatorname{ch}}
\newcommand{\Td}{\operatorname{Td}}

\newcommand{\Hilb}{\operatorname{Hilb}}
\newcommand{\Quot}{\operatorname{Quot}}

\newcommand{\Spec}{\operatorname{Spec}}
\newcommand{\Amp}{\operatorname{Amp}}

\font\b=cmr10 scaled \magstep5
\def\bigzerou{\smash{\lower1.7ex\hbox{\b 0}}}

\numberwithin{equation}{section}
\begin{document}

\title{An application of exceptional bundles to the moduli of stable sheaves
on a K3 surface}
\author{K\={o}ta Yoshioka\\
 Department of Mathematics\\
Faculty of Science, Kobe University}
 \address{Department of Mathematics, Faculty of Science,
Kobe University,
Kobe, 657, Japan }
\email{yoshioka@@math.s.kobe-u.ac.jp}
 \maketitle

\section{Introduction}
Let $X$ be a smooth projective surface defined over ${\Bbb C}$
and $L$ an ample divisor on $X$.
For a coherent sheaf $E$ on $X$,
let $v(E):=\ch(E)\sqrt{\Td(X)} \in H^*(X,{\Bbb Q})$ be 
the Mukai vector of $E$,
where $\Td(X)$ is the Todd class of $X$.
We denote the moduli of stable sheaves of
Mukai vector $v$ by $M_L(v)$,
where the stability is in the sense of Simpson [S].

For a regular surface $X$,
exceptional vector bundles
which were introduced by Drezet and Le-Potier [D-L]
 are very useful tool
for analysing $M_L(v)$.
For example, if $X={\Bbb P}^2$, 
then Maruyama [Ma1] showed the rationality
of some moduli spaces and Ellingsrud and Str\o mme [E-S]
computed relations of a generator of $H^*(M_L(v),{\Bbb Z})$.
Drezet also obtained many interesting results
[D1], [D2].
If $X$ is a K3 surface,
G\"{o}ttsche and Huybrechts [G-H] computed Hodge numbers of
rank 2 moduli spaces.
Moreover Huybrechts [H1] showed that $M_L(v)$ is deformation
equivalent to a rank 1 moduli space ( Hilbert scheme of points).
Motivated by their results,
we shall treat other rank cases. 
In particular, we shall prove the following asymptotic result.
\begin{thm}
Let $v=r+\xi+a \omega,\; \xi \in H^2(X,{\Bbb Z})$ be a primitive
Mukai vector of $r>0$ and $\langle v^2 \rangle/2 >r^2$.
Then $M_L(v)$ is deformation equivalent to 
$\Hilb_X^{\langle v^2 \rangle/2+1}$ and
$$
\theta_v:v^{\perp} \to H^2(M(v),{\Bbb Z})
$$
is an isometry which preserves hodge structures,
where $L$ is a general ample divisor.
\end{thm}
The second assertion is known by Mukai [Mu3] and
O'Grady [O] if $r \leq 2$ or $\xi$ is primitive.
For the proof of the second assertion, 
Mukai lattice and Mukai's reflection defined by
an exceptional bundle give a clear picture.

During preparation of this paper,
the author noticed that Huybrechts [H2] proved birational
irreducible symplectic manifolds are deformation equivalent.
Then primitive first Chern class cases (Theorem \ref{thm:1})
follow from O'Grady's description of $M_L(v)$ (see [H2, Cor. 4.8] and [O]).
Since our method is most successful for these cases and is needed to treat
other cases (Theorem \ref{thm:2}),
we shall first treat these cases.

\section{Preliminaries}

\subsection{Notation}

Let $M$ be a complex manifold.
For a cohomology class $x \in H^*(M,{\Bbb Z})$,
$[x]_i \in H^{2i}(X,{\Bbb Z})$ denotes the $2i$-th
component of $x$.

Let $p:X \to \Spec({\Bbb C})$ be a K3 surface over ${\Bbb C}$.
We shall recall the Mukai lattice [Mu2].
For $x,y \in H^*(X,{\Bbb Z})$, Mukai defined a symmetric 
bilinear form
$$
\langle x, y \rangle:=-p_*(x^{\vee} y) ,
$$
where $\vee:H^*(X,{\Bbb Z}) \to H^*(X,{\Bbb Z})$
be the homomorphism sending $x \in H^*(X,{\Bbb Z})
\to x-2[x]_1 \in H^*(X,{\Bbb Z})$.
For a coherent sheaf $E$ on $X$,
let $v(E):=\ch(E)(1+\omega) \in H^*(X,{\Bbb Z})$
be the Mukai vector of $E$,
where $\omega$ is the fundamental class of $X$.
Then the Riemann-Roch theorem implies that
$\chi(E,F)=-\langle v(E),v(F) \rangle$
for coherent sheaves $E$, $F$ on $X$.
Let $N$ be a line bundle on $X$.
Since $\langle x \ch(N),y \ch(N) \rangle=
\langle x,y \rangle$,
the homomorphism $T_N:H^*(X,{\Bbb Z}) \to H^*(X,{\Bbb Z})$
sending $x$ to $x \ch(N)$ is an isometry.

Let $L$ be an ample divisor on $X$.
For $v \in H^*(X,{\Bbb Z})$, 
let $M_L(v)$ be the moduli of stable sheaves
of Mukai vector $v$, where the stability is in the sense
of Simpson [S].
By Mukai [Mu1], $M_L(v)$ is smooth of dimension $
\langle v^2 \rangle+2$.
We denote the projection $S \times X \to S$ by $p_S$.
We set
$$
v^{\perp}:=\{x \in H^*(X,{\Bbb Z})| \langle v, x \rangle=0 \}.
$$
Then Mukai constructed a natural homomorphism
$$
\theta_v:v^{\perp} \longrightarrow H^2(M_L(v),{\Bbb Z})_f
$$
by
$$
\theta_v(x):=\frac{1}{\rho}
[p_{M_L(v)*}(\ch({\cal E})(1+\omega)x^{\vee})]_1,
$$
where $H^2(M_L(v),{\Bbb Z})_f$ is the torsion free quotient of
$H^2(M_L(v),{\Bbb Z})$ and ${\cal E}$ a quasi-universal family of similitude 
$\rho$.
We note that an isometry $T_N, \;N \in \Pic(X)$ satisfies that
\begin{equation}\label{eq:N}
\theta_{T_N(v)}(T_N(x))=\theta_v(x).
\end{equation}

\section{Fundamental lemma}

We shall prove the following lemma whose proof is quite
similar to [Y, Lem. 1.8]. 
\begin{lem}\label{lem:key}
Let $X$ be a smooth projective surface of $\NS(X) \cong {\Bbb Z} H$.
For a coherent sheaf $F$ of $c_1(F)=dH$, we set $\deg(F)=d$.
Let $r$ and $d$ be relatively prime positive integers
and let $r_1$ and $d_1$ be the integers which 
satisfy $r_1d-rd_1=1$ and $0<r_1 \leq r$. 
We set $r_2:=r-r_1$ and $d_2:=d-d_1$.

$(1)$
Let $E_1$ be a stable vector bundle of rank $r_1$ and $\deg(E_1)=d_1$
and $E_2$ a stable sheaf of rank $r_2$ and $\deg(E_2)=d_2$.
Then every non-trivial extension
$$
0 \to E_1 \to E \to E_2 \to 0
$$
defines a stable sheaf.

$(2)$ Let $E_1$ be a stable vector bundle of rank $r_1$ and $\deg(E_1)=d_1$
and $E$ a stable sheaf of rank $r$ and $\deg(E)=d$.
Let $V$ be a subvector space of $\Hom(E_1,E)$.
Then $V \otimes E_1 \to E$ is injective or surjective in
codimension 1.
Moreover if $V \otimes E_1 \to E$ is injective, then the
cokernel is stable.
\end{lem}
\begin{pf}
(1) We first treat the case $r_1<r$.
If $E$ is not stable, then there is a semi-stable subsheaf $G$
of $E$ such that $\deg (G)/\rk G>d/r$.
Since $G$ and $E_2$ are semi-stable and $\phi:G \to E\to E_2$ is not zero,
 $\deg (G)/\rk G \leq d_2/r_2$.
We assume that $\deg (G)/\rk G<d_2/r_2$.
Then we see that 
$1/rr_2=d_2/r_2-d/r>d_2/r_2-\deg (G)/\rk G \geq 1/r_2\rk G$,
which is a contradiction.
Hence $\deg (G)/\rk G=d_2/r_2$.
Since $G$ is semi-stable and $E_2$ is stable,
$\ker \phi$ is semi-stable of
$\deg (\ker \phi)/\rk(\ker \phi)=d_2/r_2$
and $\phi$ is surjective in codimension 1.
Hence $G$ is isomorphic to $E_2$ in codimension 1.
Let $e \in \Ext^1(E_2,E_1)$ be the extension class.
By the homomorphism $\Ext^1(E_2,E_1) \to \Ext^1(G,E_1)$,
$e$ goes to $0$.
Since $E_1$ is a vector bundle,
$\Ext^1(E_2/G,E_1)=0$.
Hence $\Ext^1(E_2,E_1) \to \Ext^1(G,E_1)$ is injective.
Thus we get that $e=0$, which is a contradiction.

If $r_1=r$, then $r=r_1=1$ and $d=d_2=1$.
In this case, it is sufficient to prove that 
$E$ is torsion free.
Let $G$ be the torsion subsheaf.
Since $E_1$ is locally free, $G \to E_2$ is injective.
Since $E_2$ is stable, $\Supp E_2/G$ is of codimension 2.
In the same way, we see that the extension class is
trivial, which is a contradiction.  

(2) We shall prove our claim by induction on $\dim V$.
We assume that $\dim V=1$.
Let $\varphi:E_1 \to E$ be a non-zero homomorphism.
We first treat the case $r_1<r$.
Since $E_1$ and $E$ are stable, $d_1/r_1 \leq 
\deg (\varphi(E_1))/\rk  \varphi(E_1)<d/r$.
In the same way as in the proof of (1),
we see that $\rk \varphi(E_1)=r_1$ and
$\deg (\varphi(E_1))=d_1$.
Hence we get that $E_1 \cong \varphi(E_1)$.
We set $E_2:=\coker \varphi$.
We assume that there is a quotient $G$ of $E_2$
such that
$G$ is semi-stable and $d_2/r_2>\deg (G)/\rk G$.
Since $G$ is a quotient of $E$,
we get that $d/r<\deg (G)/\rk G$.
Hence we get that $d/r<\deg (G)/\rk G<d_2/r_2$.
Then $1/rr_2=d_2/r_2-d/r>d_2/r_2-\deg (G)/\rk G \geq
1/r_2\rk G$,
which is a contradiction.
Hence the support of the torsion submodule of $E_2$
is 0 dimensional and the torsion free quotient is stable.
Since $E_1$ is locally free, $E_2$ is torsion free and
hence stable.
If $r_1=r \;( \;= 1\;)$, then $\varphi$ is injective.
We shall show that $\coker \varphi$ is of pure dimension $1$.
Let $T$ be a subsheaf of dimension 0 and $\widetilde{T}$
the pull-back of $T$ to $E$.
Then $E_1 \to \widetilde{T}$ is an isomorphism
in codimension 1.
Since $E_1^{\vee \vee} \cong \widetilde{T}^{\vee \vee}$
and $E_1$ is locally free, $E_1=\widetilde{T}$.
Hence we get that $T=0$.
We set $C:=\Supp E_2$.
Since $H$ is a generator of the N\'{e}ron-Severi group,
$C$ is reduced and irreducible.
Hence $E_2$ is a torsion free sheaf of rank 1 on $C$,
which shows that $E_2$ is a stable sheaf.

We shall treat general cases.
Let $\varphi_1,\varphi_2,\dots,\varphi_n$ be a basis of $V$
and $V'$ the subspace generated by $\varphi_2,\dots,\varphi_n$.
We note that $r_1d_2-r_2d_1=1$.
We assume that $r_1<r_2$.
By induction hypothesis,
$V' \otimes E_1 \to \coker(\varphi_1)$ is injective or
surjective in codimension 1.
If $r_1 > r_2>0$, in the same way, we see that $\varphi_2:E_1 \to 
\coker(\varphi_1)$ is surjective.
If $r_1=r_2$ or $r_2=1$, then $r_1=d_2=1$.
Since every degree 1 curve is reduced and irreducible,
 $V' \otimes E_1 \to \coker(\varphi_1)$ is injective or
surjective in codimension 1. 
Hence we get our claim. 
\end{pf}

\begin{cor}\label{cor:uni}
Under the same assumption of Lemma \ref{lem:key} $(1)$,
the universal extension
$$
0 \to E_1\otimes \Ext^1(E_2,E_1)^{\vee} \to E
\to E_2 \to 0
$$
defines a stable sheaf.
\end{cor}

\section{Moduli of stable sheaves on K3 surfaces}

\subsection{Stable sheaves of pure dimension 1}
Let $H$ be an ample divisor on $X$.
We set $v:=H+a \omega \in H^*(X,{\Bbb Z})$.
For a pure dimensional sheaf $E$ of $v(E)=v$ and
an ample divisor $L$,
$\chi(E \otimes L^{\otimes n})=(H,L)n+b$.
Hence $E$ is semi-stable with respect to $L$, if
$$
\frac{\chi(F)}{(c_1(F),L)} \leq \frac{\chi(E)}{(c_1(E),L)} 
$$
for any proper subsheaf $F$ of $E$,
and $E$ is stable if the inequality is strict.
We can generalize the concept of the chamber.
For a sheaf $F$ of pure dimension 1,
We set $D:=\chi(F)c_1(E)-\chi(E)c_1(F)$ and
$W_D:=\{x \in \Amp(X)|(x,D)=0\}$.
Then $(D^2)=\chi(F)^2(c_1(E)^2)-2\chi(E)\chi(F)(c_1(E),c_1(F))
+\chi(E)^2(c_1(F)^2)$.
If $W_D$ is not empty,
then the Hodge index theorem implies that $(D^2) \leq 0$.
Hence the choice of $\chi(F)$ is finite,
which shows that the number of
non-empty walls $W_D$ is finite.
We shall call a chamber a connected component of $\Amp(X) \setminus
\cup_D W_D$.
If $c_1(E)$ is primitive, 
then $M_L(v)$ is compact for a general $L$.
As is well known ([G-H, Prop. 2.2], cf. [Y, Prop. 3.3]), the choice of polarizations 
is not so important.
Hence we denote $M_L(v)$ by $M(v)$.

\subsection{Correspondence}

Let $r$ and $d$ be relatively prime positive integers
and let $r_1$ and $d_1$ be the integers which 
satisfy $r_1d-rd_1=1$ and $0<r_1 \leq r$.
We shall consider a K3 surface $X$ such that
$\Pic(X)= {\Bbb Z} H$, where $H$ is an ample divisor.
We assume that there are Mukai vectors
$v_1,v \in H^*(X,{\Bbb Z})$ such that
\begin{equation}
\begin{cases}
v_1=r_1+d_1H+a_1 \omega,\\
v=r+dH+a \omega,\\
\langle v_1^2 \rangle =-2,
\end{cases}
\end{equation}
where $a_1,a \in {\Bbb Z}$.
Since $\langle v_1^2 \rangle=-2$,
there is a unique stable vector bundle $E_1$
of $v(E_1)=v_1$, that is $M(v_1)=\{E_1 \}$.
We note that $E_1$ satisfies that
\begin{equation}
\begin{cases}
\Hom(E_1,E_1)={\Bbb C},\\
\Ext^1(E_1,E_1)=0,\\
\Ext^2(E_1,E_1)={\Bbb C}.
\end{cases}
\end{equation}

\begin{defn}
For $i \geq 1$,
$$
M(v)_i:=\{E \in M(v)| \dim \Hom(E_1,E)=-\langle v_1,v \rangle-1+i \}
$$
is a locally closed subscheme of $M(v)$ with reduced structure.
\end{defn}

\begin{defn}
For $w:=v-mv_1$ with $[w]_0 \geq 0$,
we set
$$
N(mv_1,v,w):=\{E_1^{\oplus m} \subset E| E \in M(v) \}.
$$
Let $\pi_v:N(mv_1,v,w) \to M(v)$ is the projection and 
$N(mv_1,v,w)_i=\pi_v^{-1}(M(v)_i)$.
\end{defn}
We shall construct $N(mv_1,v,w)$ as a scheme over $M(v)$.
Let $Q(v)$ be the open subscheme of a canonical quot scheme
$\Quot_{{\cal O}_X^{\oplus N}/X}$ such that 
$Q(v)/PGL(N) \cong M(v)$
and ${\cal O}_{Q(v) \times X}^{\oplus N} \to
 {\cal E}_v$ be the universal quotient on $Q(v) \times X$.
Let $W_2 \to W_1 \to E_1 \to 0$ be a locally free resolution
of $E_1$ such that
$\Ext_{p_{Q(v)}}^i(W_j \boxtimes {\cal O}_{Q(v)},{\cal E}_v)=0$ 
for $i>0$ and $j=1,2$.
We set ${\cal W}_j:=
\Hom_{p_{Q(v)}}(W_j \boxtimes {\cal O}_{Q(v)},{\cal E}_v)$, $j=1,2$.
We shall consider the Grassmannian bundle
$\xi:{\Bbb G}=Gr({\cal W}_1,m) \to Q(v)$ parametrizing $m$-dimensional
subspaces of $({\cal W}_1)_t, t \in Q(v)$.
Let ${\cal U} \to \xi^*{\cal W}_1$ be the universal subbundle.
For the composition
 $\Psi:{\cal U} \to \xi^*{\cal W}_1 \to \xi^*{\cal W}_2$,
we set $D:=\{t \in {\Bbb G}| \Psi_t=0\}$.
Then there is a universal family of homomorphisms on $D \times X$:
$$
E_1 \boxtimes {\cal U} \to \xi^*{\cal E}_v.
$$  
By Lemma \ref{lem:key},
$m$-dimensional subspace $U$ of
$\Hom(E_1,E)$ defines an injective homomorphism 
$E_1 \otimes U \to E$.
Since ${\cal W}_1,{\cal W}_2, {\cal U}$ and $\Psi$ are
$GL(N)$-equivariant, We get $N(mv_1,v,w) \to M(v)$
as a scheme $D/PGL(N) \to M(v)$. 
Since $F:=\coker(E_1 \otimes U \to E)$ is a stable sheaf,
we also obtain a morphism
$\pi_w:N(mv_1,v,w) \to M(w)$.
We shall consider the fiber of $\pi_w$.
For $F \in M(w)_{i+m}$, we get that
$\dim \Ext^1(F,E_1)=i-1+m$.
Let $U^{\vee}$ be a $m$-dimensional subspace
of $\Ext^1(F,E_1)$.
Then the element $e \in U \otimes \Ext^1(F,E_1)$
which corresponds to the inclusion $U^{\vee} \subset \Ext^1(F,E_1)$
defines an extension
$$
0 \to U \otimes E_1 \to E \to F \to 0.
$$
By the choice of the extension class and Corollary \ref{cor:uni},
$E$ is a stable sheaf.
Thus the fiber of $F$ is the Grassmannian $Gr(i-1+m,m)$.
 
\begin{lem}\label{lem:cod}
For $i \geq 1+\langle v,v_1 \rangle$,
$$
\codim M(v)_i=(i-1)(i-1-\langle v,v_1 \rangle).
$$
In particular, if $\langle v,v_1 \rangle \leq 0$, then 
$M(v)_1$ is an open and dense subscheme of
$M(v)$.
\end{lem}

\begin{pf}
Let $E$ be an element of $M(v)_i$.
Then $\dim \Ext^1(E_1,E)=-\chi(E_1,E)+\dim \Hom(E_1,E)
=i-1$.
We shall consider the universal extension
$$
0 \to E_1^{\oplus (i-1)} \to G \to E \to 0.
$$
By Corollary \ref{cor:uni},
$G$ is a stable sheaf of $\rk(G)=(i-1)r_1+r$ and
$\deg(G)=(i-1)d_1+d$.
Hence we get that $\Ext^2(E_1,G)=0$.
Simple calculations show that
\begin{align*}
\dim \Hom (E_1,G)&=i-1+\dim \Hom(E_1,E)\\
&=2i-2-\langle v_1,v \rangle,\\
\dim \Ext^1(E_1,G)&=0.
\end{align*}
Hence $G$ belongs to $M(v(G))_1$.
We set $k:=2i-2-\langle v_1,v \rangle$.
Then the choice of $i-1$-dimensional subspace
of $V:=\Hom (E_1,G)$ is parametrized by
the Grassmannian $Gr(k,i-1)$.
Since $v(G)=v+(i-1)v_1$,
we see that
\begin{align*}
\dim M(v(G))&=\langle v(G)^2 \rangle+2\\
&=-2(i-1)^2+2(i-1)\langle v_1,v \rangle +\langle v^2 \rangle+2.
\end{align*}
Therefore we see that
\begin{align*}
\dim M(v)_i &=\dim M(v(G))+\dim Gr(k,i-1)\\
&=-(i-1)(i-1-\langle v_1,v \rangle)+\dim M(v),
\end{align*}
which implies our claim.
\end{pf}

\begin{rem}\label{rem:gr}
We can also see that $M(v)_i$ is an \'{e}tale locally trivial
$Gr(k,i-1)$-bundle over $M(v(G))_1$.
Moreover we can show that
$\pi_v:N(mv_1,v,w)_i \to M(v)_i$ is an \'{e}tale 
locally trivial $Gr(l,m)$-bundle and
$\pi_w:N(mv_1,v,w)_i \to M(w)_{i+m}$ is an \'{e}tale 
locally trivial $Gr(i-1+m,m)$-bundle,
where $l:=i-1-\langle v_1,v \rangle$.

\end{rem}

\subsection{A special case}

We set $k(s):=s r_1^2+rr_1-r^2$, $s \in {\Bbb Z}$.
For a positive integer $k(s)$, we shall consider a K3 surface
$X$ such that $\Pic(X)={\Bbb Z} H$ of $(H^2)=2k(s)$.
We set
\begin{equation}
\begin{cases}
v_1=r_1+d_1H+(d_1^2r+d_1^2sr_1-r_1d^2+2d)\omega\\
v=r+dH+((2dd_1r_1-rd_1^2)s+d^2(r_1-r))\omega.
\end{cases}
\end{equation}
Then a simple calculation shows that
\begin{equation}
\begin{cases}
\langle v_1^2 \rangle=-2\\
\langle v^2 \rangle=2s\\
\langle v,v_1 \rangle=-1.
\end{cases}
\end{equation}

\begin{lem}\label{lem:1}
We assume that $r_1 \geq r/2$.
Under the following assumptions, $k(s)>0$.
\begin{enumerate}
\item
$r=2$ and $s \geq 3$,
\item
$r_1 \geq (2/3) r$ and $s \geq 1$.
\end{enumerate}
\end{lem}

Let us consider the reflection $R_{v_1}$ of
$H^*(X,{\Bbb Z})$ defined by $v_1$:
$$
x \mapsto x+\langle x,v_1 \rangle v_1.
$$
We set $w:=R_{v_1}(v)=v-v_1$.
By Lemma \ref{lem:cod}, $N(v_1,v,w)$ gives a birational correspondence
between $M(v)$ and $M(w)$.

\begin{lem}\label{lem:smoo}
$N(v_1,v,w)$ is smooth and irreducible.
\end{lem}
\begin{pf}
In the same way as in [G-H],
it is sufficient to prove that
$\Hom(E,E_2)={\Bbb C}$.
We assume that $\dim \Hom(E,E_2)>1$.
Let $\phi:E \to E_2$ be a non-zero
homomorphism.
By Lemma \ref{lem:key}, $\phi$ is surjective in codimension 1
and we also see that $\ker \phi$ is stable.
Since $[v(\ker \phi)]_i=[v_1]_i$, $i=0,1$, the stability of $\ker \phi$
implies that  
$[v(\ker \phi)]_2 \leq [v_1]_2$.
Since $[v(\phi(E))]_2 \leq [w]_2$ and $[v(\phi(E))]_2+[v(\ker \phi)]_2
=[v]_2$, we must have
$v(\ker \phi)=v_1$ and $v(\phi(E)) =w$.
Hence $\phi$ is surjective and $\ker \phi \cong E_1$.
Let 
$$
0 \to E_1 \otimes \Ext^1(E,E_1)^{\vee} \to G \to E \to 0
$$
be the universal extension.
Then $\ker \phi$ is determined by a $m$-dimensional 
subspace of $ \Hom(E_1,G)$ containing $\Ext^1(E,E_1)^{\vee}$,
where $m=\dim \Ext^1(E,E_1)+1$.
We shall consider the composition $\phi':G \to E \to E_2$.
Then it defines a universal extension (up to the action of
$\Aut(\Ext^1(E_2,E_1)^{\vee})$)
$$
0 \to E_1 \otimes \Ext^1(E_2,E_1)^{\vee} \to G \to E_2 \to 0.
$$
Since $G$ is stable,
$\phi'$ is determined by a subspace $\Ext^1(E_2,E_1)^{\vee}$
of $\Hom(E_1,G)$ up to multiplication by constants. 
Hence $\dim \Hom(E,E_2)=1$.
\end{pf}

Let us consider the relation
between $\theta_v$ and $\theta_{w}$.
By our assumption on $v$, there is a universal
family ${\cal E}$ on $M(v) \times X$ ([M1]).
Then $N(v_1,v,w)$ is constructed as a projective bundle
${\Bbb P}(V)$, where 
$V:=\Ext^2_{p_{M(v)}}({\cal E},{\cal O}_{M(v)}\boxtimes(E_1 \otimes K_X))$,
and there is a universal homomorphism 
$\Psi:{\cal O}_{N(v_1,v,w)}\boxtimes E_1 \to
\pi_v^*{\cal E} \otimes {\cal O}_{N(v_1,v,w)}(1)$. 
 We denote the cokernel by ${\cal G}$.
Then ${\cal G}$ is flat over $N(v_1,v,w)$.
Since $\codim(M(v)_i) \geq 2$
for $i \geq 2$,
we see that $\det p_{N(v_1,v,w)!}
({\cal G} \otimes q^*(E_1^{\vee})) \cong {\cal O}_{N(v_1,v,w)}$.
Therefore we see that
\begin{align*}
[p_{N(v_1,v,w)*}(\ch({\cal E}\otimes {\cal O}_{N(v_1,v,w)}(1))
(1+\omega)x^{\vee})]_1
& =[p_{N(v_1,v,w)*}(\ch({\cal G})(1+\omega)x^{\vee})]_1\\
& =[p_{N(v_1,v,w)*}(\ch({\cal G})(1+\omega)
(x+\langle v_1,x \rangle v_1)^{\vee})]_1.
\end{align*}
Hence we get that 
$\theta_v(x)=\theta_{R_{v_1}(v)}(R_{v_1}(x))$
for $x \in v^{\perp}$.

\begin{prop}\label{prop:bilin}
The following diagram is commutative.
\begin{equation*}
\begin{CD}
v^{\perp} @>{R_{v_1}}>> w^{\perp}\\
@V{\theta_v}VV  @VV{\theta_{w}}V\\
H^2(M(v),{\Bbb Z})_f @>>{\pi_{w*}\pi_v^*}>
H^2(M(w),{\Bbb Z})_f
\end{CD}
\end{equation*}
\end{prop}

Let us consider the structure of
$N(v_1,v,w)_i$ more closely.
For simplicity we assume that $i=2$.
We set $u:=v+v_1$.
Let $Q(u)$ be the open subscheme of a quot scheme
$\Quot_{{\cal O}_X^{\oplus N}/X}$ such that $Q(u)/PGL(N) \cong M(u)$
and let ${\cal O}_{Q(u) \times X}^{\oplus N} \to
 {\cal E}_u$ be the universal quotient on $Q(u) \times X$.
We set $V:=\Hom_{p_{Q(u)_1}}(E_1 \boxtimes {\cal O}_{Q(u)_1},{\cal E}_u)$,
where $Q(u)_1$ is the pull-back of $M(u)_1$ to $Q(u)$.
Then $V$ is a locally free sheaf of rank $3$.
We shall consider projective bundles
$\xi_1:{\Bbb P}_1:={\Bbb P}(V^{\vee}) \to Q(u)_1$and 
$\xi_2:{\Bbb P}_2:={\Bbb P}(V)\to Q(u)_1$.
Let
\begin{align*}
&0 \to {\cal Q}_1 \to \xi_1^*V^{\vee} \to {\cal O}_{{\Bbb P}_1}(1) \to 0\\
&0 \to {\cal Q}_2 \to \xi_2^*V \to {\cal O}_{{\Bbb P}_2}(1) \to 0
\end{align*}
be the universal bundles on ${\Bbb P}_1$ and ${\Bbb P}_2$ 
respectively.
Let ${\cal Z} \subset {\Bbb P}_1 \times {\Bbb P}_2$ be the 
incidence correspondence.
Let $\eta_i:{\cal Z} \to {\Bbb P}_i$ $i=1,2$ be the projections. 
For simplicity, we denote $\eta_1^* {\cal O}_{{\Bbb P}_1}(m) \otimes
\eta_2^* {\cal O}_{{\Bbb P}_2}(n)$ by ${\cal O}_{{\cal Z}}(m,n)$.
There is a universal family of filtrations
\begin{equation}
0 \subset {\cal O}_{{\cal Z}}(-1,0) \subset \eta^*_2 {\cal Q}_2
\subset \eta_2^* \xi_2^* V.
\end{equation}
Then there are homomorphisms 
\begin{align*}
\alpha:& E_1 \boxtimes {\cal O}_{{\Bbb P}_1}
\to E_1 \boxtimes \xi_1^*V \otimes {\cal O}_{{\Bbb P}_1}(1)
\to (\xi_1 \times id_X)^*{\cal E}_u \otimes
{\cal O}_{{\Bbb P}_1 \times X}(1) \\
\beta:& E_1 \boxtimes {\cal Q}_2
\to E_1 \boxtimes \xi_2^*V
\to (\xi_2 \times id_X)^*{\cal E}_u.
\end{align*}
By the construction of the homomorphisms,
$\alpha$ and $\beta$ are injective and ${\cal E}_v:=\coker \alpha$
and ${\cal E}_w:=\coker \beta$ are flat family of stable sheaves of 
Mukai vectors $v$ and $w$ respectively.
Then we get the following exact and commutative diagram.
\begin{equation}
\begin{CD}
@. @. 0 @. 0 @.\\
@.@. @VVV @VVV\\
0 @>>> E_1 \boxtimes {\cal O}_{{\cal Z}}(-1,0) @>>>
E_1\boxtimes \eta^*_2 {\cal Q}_2 @>>>
E_1 \boxtimes {\cal O}_{{\cal Z}}(1,-1) @>>> 0 \\
@. @| @VVV @VV{\psi}V @. \\
0 @>>> E_1 \boxtimes {\cal O}_{{\cal Z}}(-1,0) @>>>
(\xi_1 \eta_1 \times id_X)^*{\cal E}_u
@>>> (\eta_1 \times id_X)^*{\cal E}_v @>>> 0\\
@.@. @VVV @VVV @.\\
@.@. (\eta_2 \times id_X)^*{\cal E}_w @= 
(\eta_2 \times id_X)^*{\cal E}_w @.\\
@.@.@VVV @VVV@.\\
@. @. 0 @. 0 @.
\end{CD}
\end{equation}
The homomorphism $\psi$ defines a morphism ${\cal Z} \to N(v_1,v,w)_2$.
Obviously this morphism is $PGL(N)$ invariant,
and hence we get a morphism $f:{\cal Z}/PGL(N) \to N(v_1,v,w)$.
Conversely, let ${\cal F}_v$ be a family of stable sheaves 
which belong to $M(v)_2$ and
let $\psi:E_1 \boxtimes {\cal O}_S \to {\cal F}_v$ 
be a family of homomorphisms
which belong to $N(v_1,v,w)_2$ and are parametrized by a scheme $S$
over $M(v)_2$.
Then $L:=\Ext^1_{p_S}({\cal F}_v,E_1 \boxtimes {\cal O}_S)$
is a line bundle. 
Hence we get an extension
$$
0 \to E_1 \boxtimes L^{\vee} \to
{\cal F}_u \to {\cal F}_v \to 0.
$$
We set ${\cal F}_w:=\coker \psi$ 
and ${\cal G}:=\ker({\cal F}_u \to {\cal F}_w)$.
Then ${\cal F}_w$ is a family of stable sheaves
which belong to $M(w)$.
We get the following exact and commutative diagram.

\begin{equation}
\begin{CD}
@. @. 0 @. 0 @.\\
@.@. @VVV @VVV\\
0 @>>> E_1 \boxtimes L^{\vee} @>>>
 {\cal G} @>>>
E_1 \boxtimes {\cal O}_{S} @>>> 0 \\
@. @| @VVV @VV{\psi}V @. \\
0 @>>> E_1 \boxtimes  L^{\vee} @>>>
{\cal F}_u
@>>> {\cal F}_v @>>> 0\\
@.@. @VVV @VVV @.\\
@.@. {\cal F}_w @= 
{\cal F}_w @.\\
@.@.@VVV @VVV@.\\
@. @. 0 @. 0 @.
\end{CD}
\end{equation}
We set $Q:=\Hom_{p_S}(E_1 \boxtimes {\cal O}_S,{\cal G})$.
Then it is easy to see that $E_1 \boxtimes Q \to {\cal G}$
is an isomorphism.
$E_1 \boxtimes L^{\vee} \to {\cal F}_3$ defines a morphism 
$S \to N(v_1,u,v)_1$.
The filtration of vector bundles
$$
0 \subset L^{\vee} \subset Q \subset 
\Hom_{p_S}(E_1 \boxtimes {\cal O}_S,{\cal F}_u)
$$
defines a lifting $S \to {\cal Z}/PGL(N)$.
In particular, we obtain a morphism
$g:N(v_1,v,w)_2 \to {\cal Z}/PGL(N)$.
Then $g$ is the inverse of $f$.
We also see that ${\Bbb P}_1/PGL(N)=N(v_1,u,v)_1 \to M(v)_2$
and ${\Bbb P}_2/PGL(N)=N(2v_1,u,w)_1 \to M(w)_3$ are 
isomorphisms.
Therefore $\pi_{w}\pi^{-1}_v:M(v)\; - - > M(w)$
is an elementary transformation in codimension 2.
By using [H1, Cor. 5.5] or [H2, Cor. 4.7], we obtain the following theorem. 
\begin{thm}\label{thm:corres}
$M(v)$ is deformation equivalent to
$M(w)$.
\end{thm}

\begin{rem}
In the notation of Lemma \ref{lem:cod} and Remark \ref{rem:gr},
we shall consider \'{e}tale locally trivial 
$Gr(2i-1,i-1) \times Gr(2i-1,i)$-bundle 
$M(v)_i \times_{M(v(G))_1} M(w)_{i+1} \to M(v(G))_1$.
Let ${\cal Z}$ be the incidence correspondence. Then we also see that
${\cal Z}$ is isomorphic to $N(v_1,v,w)_i$.
\end{rem}

\subsection{Cohomologies of $M(v)$}

\begin{lem}\label{lem:trans}
Assume that $\rho(X) \geq 2$.
Let $v=l(r+\xi)+a \omega,\; \xi \in H^2(X,{\Bbb Z})$ be a Mukai vector
such that $r+ \xi$ is primitive.
Then there is a line bundle $L$ such that
$(1)$ $\xi':=[T_L(v)]_1/l=\xi+rL$ is primitive, 
$(2)$ $\xi'$ is ample, and
$(3)$ $({\xi'}^2) \geq 4$.
\end{lem}
\begin{pf}
Let $L'$ be a primitive ample divisor on $X$ such that
$L'$ and $\xi$ are linearly independent.
Let $n$ be an integer
such that $r$ and $n$ are relatively prime.
Since $L'$ is primitive, and
$[T_{nL'}(v)]_1/l=\xi+rnL'$,
$[T_{nL'}(v)]_1/l$ is primitive.
Hence for a sufficiently large integer $n$,
$L:=nL'$ satisfies our claims.
\end{pf}

\begin{prop}\label{prop:deform}
Let $X_1$ and $X_2$ be K3 surfaces, and let
$v_1:=l(r+\xi_1)+a_1 \omega \in H^*(X_1,{\Bbb Z})$
and $v_2:=l(r+\xi_2)+a_2 \omega \in H^*(X_2,{\Bbb Z})$
be Mukai vectors such that
$(1)$ $r+\xi_1$ and $r+\xi_2$ are primitive,
$(2)$ $\langle v_1^2 \rangle=\langle v_2^2 \rangle=2s$,
and $(3)$ $a_1 \equiv a_2 \mod l$.
Then $M(v_1)$ and $M(v_2)$ are deformation equivalent.

\end{prop}

\begin{pf}
We may assume that $\xi_1$ and $\xi_2$ are ample.
We assume that $\rho(X_i)=1$ for some $i$.
Let $T_i$ be a connected smooth curve and
$({\cal X}_i, {\cal L}_i)$ a pair of a smooth family of
K3 surfaces $p_{T_i}:{\cal X}_i \to T_i$
 and a relatively ample line bundle
${\cal L}_i$.
For points $t_0, t_1 \in T_i$,
 we assume that $(({\cal X}_i)_{t_0},({\cal L}_i)_{t_0})
=(X_i,\xi_i)$ and $({\cal X}_i)_{t_1}$ is a K3 surface of 
$\rho(({\cal X}_i)_{t_1}) \geq 2$.
We can construct an algebraic space 
${\cal M}_{{\cal X}_i/T_i}(v_i) \to T_i$
which is smooth and proper over $T_i$, and is a family
of moduli of stable sheaves on geometric fibers of Mukai vector $v_i$ 
with respect to general polarizations on fibers ([G-H],[O],[Y]).
Replacing $X_i$ by $({\cal X}_i)_{t_1}$, we may assume that
$\rho(X_i) \geq 2$ for $i=1,2$.
By Lemma \ref{lem:trans}, we may assume that 
$(1)$ $\xi_i$ is primitive, 
$(2)$ $\xi_i$ is ample, and
$(3)$ $(\xi_i^2) \geq 4$.
Let $\pi:Z \to {\Bbb P}^1$ be an elliptic K3 surface of $\rho(Z)=2$ 
and let $f$ be a fiber of $\pi$ and $C$ a section of $\pi$.
Then by using deformations of $(X_i,\xi_i)$ and $M(v_i)$, 
we see that $M(v_i)$ is deformation equivalent to
$M(v_i')$, where 
$v_i':=l(r+(C+k_i f))+a_i \omega$, $i=1,2$.
Since $a_1 \equiv a_2 \mod l$,
we set $n:=(a_1-a_2)/l$.
Then it is easy to see that $T_{n f}(v_2')=v_1'$.
\end{pf}

\begin{rem}\label{rem:theta}
By the proof of [Mu2, Thm. A.5], 
there is a quasi-universal family ${\cal F}$ on
${\cal M}_{{\cal X}_i/T_i}(v_i) \times_{T_i} {\cal X}_i$.
Since ${\cal L}$ is relatively ample, we can consider 
a locally free resolution of ${\cal F}$.
In particular, we can construct a family of homomorphisms
$(\theta_{v_i})_t:(v_i)^{\perp}_t \to 
H^2({\cal M}_{({\cal X}_i)_t/k(t)}(v_i),{\Bbb Z})$,
where $(v_i)^{\perp}_t \subset H^*(({\cal X}_i)_t,{\Bbb Z})$.
This fact will be used in the proof of Corollary 
\ref{cor:period} (cf. [Y, Prop. 3.3]).
\end{rem}

We get another proof of [H2, Cor. 4.8].
\begin{thm}\label{thm:1}
Let $v=r+ \xi+a \omega,\; \xi \in H^2(X,{\Bbb Z})$ be a Mukai vector
such that $r>0$ and $r+\xi$ is primitive.
Then $M(v)$ and 
$\Hilb_X^{\langle v^2 \rangle/2+1}$ are deformation equivalent.
If $r=0$, $\xi$ is ample, and $(\xi^2) \geq 4$,
 then the same results hold.  
\end{thm}

\begin{pf}
If $\langle v^2 \rangle =0$, then Mukai [Mu1] showed that 
$M(v)$ is a K3 surface.
Hence we assume that $\langle v^2 \rangle =2s>0$.
We first assume that $r \ne 2$.
Let $v$ be a Mukai vector in Theorem \ref{thm:corres}.
We assume that $d=r-1$. Since $r_1=r-1$,
$M(v)$ is deformation equivalent to $\Hilb_X^{s+1}$.
By Proposition \ref{prop:deform},
$M(v)$ is deformation equivalent to $\Hilb_X^{s+1}$ for every $v$.

We shall treat the rank 2 case.
Let $s$ be a positive integer.
By Theorem \ref{thm:corres} and Lemma \ref{lem:1},
there is a Mukai vector $v$ and $w$ such that
$\langle v^2 \rangle=\langle w^2 \rangle=2s$,
$[v]_0=2$, $[w]_0=7$, and $M(v)$ is deformation equivalent
to $M(w)$.
By using Proposition \ref{prop:deform},
 we see that $M(v)$ is
deformation equivalent to
$\Hilb_X^{\langle v^2 \rangle/2+1}$
foe every $v$.
\end{pf}

We also get another proof of [O].
\begin{cor}\label{cor:period}
Under the same assumption of Theorem \ref{thm:1} 
and $\langle v^2\rangle \geq 2$,
$M(v)$ is an irreducible symplectic manifold and
$$
\theta_v : v^{\perp} \longrightarrow H^2(M(v),{\Bbb Z})
$$
is an isometry which preserves hodge structures.
\end{cor}
\begin{pf}
The rank one case easily follows from [B]
(cf. [Mu3],[O]).
By using \eqref{eq:N}, Proposition \ref{prop:bilin},
Remark \ref{rem:theta},
and the proofs of Proposition \ref{prop:deform} and 
Theorem \ref{thm:1},
we get our corollary.
\end{pf}

\section{Non-primitive first Chern class cases}
In this section, we shall consider a more general case 
and get some partial results.
Let $r$ and $d$ be relatively prime non-negative integers.
Then there are integers $r_1$ and $d_1$ such that
$dr_1-d_1r=1$
Let $l$ be a positive integer such that 
$l$ and $r_1$ are relatively prime.
We shall choose $r_1$ of $0<r_1<lr$.
Then there are unique pair of integers $p,q$ of
$pr_1-ql=-1$ and $0 \leq p < l$.
We set $k(s):=r_1(qr+r_1s)-r^2, s \in {\Bbb Z}$.
Let $X$ be a K3 surface such that 
$\Pic(X)={\Bbb Z} H$ and $(H^2)=2k(s)$.
We set
\begin{equation}\label{eq:4-1}
\begin{cases}
v=lr+ld H+\{l((1+dr_1)d_1s+d^2qr_1-rd^2)-p \}\omega,\\
v_1=r_1+d_1 H+\{r_1(-d^2+d_1^2s)+d_1^2rq+2d \}\omega.
\end{cases}
\end{equation}
Then we see that
\begin{equation}
\begin{cases}
\langle v_1^2 \rangle=-2,\\
\langle v^2 \rangle =2l(ls+rp),\\
\langle v_1,v \rangle=-1.
\end{cases}
\end{equation}
We set $w:=v-v_1$ and we shall consider the relation between
$M(v)$ and $M(w)$.

\begin{lem}\label{lem:a1}
Let $E_1$ be a stable vector bundle of $\rk(E_1)=r_1$ and
$\deg(E_1)=d_1$.

$(1)$ Let $E$ be a $\mu$-stable sheaf of $\rk(E)=lr$ and 
$\deg(E)=ld$.
Then every non-zero homomorphism $\phi:E_1 \to E$ is 
injective and $\coker \phi$ is a stable sheaf.

$(2)$ Let $E$ be a $\mu$-stable sheaf of $\rk(E)=lr$ and 
$\deg(E)=ld$.
For a non-trivial extension
$$
0 \to E_1 \to E' \to E \to 0,
$$
$E'$ is a stable sheaf.

$(3)$ Let $E_2$ be a stable sheaf of 
$\rk(E_2)=lr-r_1$ and $\deg(E_2)=ld-d_1$.
Then every non-trivial extension
$$
0 \to E_1 \to E \to E_2 \to 0,
$$
defines a $\mu$-semi-stable sheaf.

$(4)$ Let $E'$ be a stable sheaf of $\rk(E')=lr+r_1$
and $\deg(E')=ld+d_1$.
Then every non-zero homomorphism $\phi:E_1 \to E'$ is 
injective and $\coker \phi$ is $\mu$-semi-stable.
\end{lem}

\begin{pf}
(1) Since $d_1/r_1=\deg(E_1)/\rk(E_1) \leq 
\deg(\phi(E_1))/\rk(\phi(E_1))
 < \deg(E)/\rk(E)=d/r$,
we get that 
$1/rr_1=d/r-d_1/r_1 \geq d/r-\deg(\phi(E_1))/\rk(\phi(E_1))
 \geq 1/r \rk(\phi(E_1))$.
Hence we obtain that $\rk(\phi(E_1))=r_1$, which implies that
$\phi$ is injective.
We assume that $\coker \phi$ is not stable.
Then there is a semi-stable quotient sheaf $G$ of 
$\coker \phi$ such that
$\deg(G)/\rk(G) <\deg(\coker \phi)/\rk(\coker \phi)
=(ld-d_1)/(lr-r_1)$.
Since $G$ is a quotient of $E$, we get that
$\deg(G)/\rk(G)>d/r$.
Hence we obtain that $1/(lr-r_1)r=
(ld-d_1)/(lr-r_1)-d/r>\deg(G)/\rk(G)-d/r 
\geq 1/r \rk(G)$, which is a contradiction.
Therefore $G$ is a stable sheaf.

(2) We assume that $E'$ is not stable.
Let $G$ be a destabilizing semi-stable subsheaf of $E'$.
We assume that $\phi:G \to E$ is not surjective in codimension 1.
Then the $\mu$-stability of $E$ implies that
 $d/r > \deg(\phi(G))/\rk(\phi(G)) \geq 
\deg(G)/\rk(G)>\deg(E')/\rk(E')=(ld+d_1)/(lr+r_1)$.
Hence we see that
$1/r(lr+r_1)>d/r-\deg(G)/\rk(G) 
\geq 1/r \rk(G)$, which is a contradiction.
Thus $\phi$ is surjective in codimension 1.
If $\deg(\phi(G))/\rk(\phi(G)) >
\deg(G)/\rk(G)$, then we also get a contradiction,
and hence $\ker \phi$ is $\mu$-semi-stable of 
$\deg(\ker \phi)/\rk(\ker \phi)=\deg(G)/\rk(G)$.
Thus $\deg(\ker \phi)/\rk(\ker \phi)>d_1/r_1$,
which contradicts the stability of $E_1$.
Therefore $G \to E$ is an isomorphism in codimension 1.
In the same way as in the proof of
Lemma \ref{lem:key}, we see that the extension is split,
which is a contradiction.

(3) We assume that $E$ is not $\mu$-semi-stable.
Let $G$ be a semi-stable subsheaf of $E$ such that
$d/r<\deg(G)/\rk(G)$.
It is sufficient to show that $\phi:G \to E_2$
is an isomorphism in codimension 1.
We note that
$d/r<\deg(G)/\rk(G) \leq \deg(\phi(G))/\rk(\phi(G))
\leq (ld-d_1)/(lr-r_1)$.
Then we get that 
$1/r(lr-r_1)=(ld-d_1)/(lr-r_1)-d/r \geq 
\deg(\phi(G))/\rk(\phi(G))-d/r \geq 1/\rk(\phi(G))r$.
Thus we obtain that $\rk(\phi(G))=(lr-r_1)$ and 
$\deg(\phi(G))/\rk(\phi(G))
= (ld-d_1)/(lr-r_1)$.
Therefore $\phi$ is surjective in codimension 1.
We set $G':=\ker \phi$.
We assume that $G' \ne 0$.
Since 
$\rk(G)=\rk(G')+lr-r_1$,
$\deg(G)=\deg(G')+ld-d_1$ and
$d/r-\deg(G)/\rk(G)<0$,
we get that $r \deg(G') \geq d \rk(G')$.
Hence $\deg(G')/\rk(G') > \deg(E_1)/\rk(E_1)$,
which is a contradiction.
Therefore $G \to E_2$ is an isomorphism in codimension 1.

The proof of (4) is similar to that
of (1). 
\end{pf}

Let $M(v)^{\mu}$ be the open subscheme of $M(v)$
consisting of $\mu$-stable sheaves. 
In the same way as in section 3, we shall define 
$N(v_1,v,w)$ and an open subscheme 
$N(v_1,v,w)^{\mu}:=\{E_1 \subset E| E \in M(v)^{\mu} \}$.

\begin{lem}\label{lem:smooth2}
$N(v_1,v,w)^{\mu}$ is smooth and irreducible.
\end{lem}
\begin{pf}
For $E_1 \subset E \in N(v_1,v,w)^{\mu}$,
we set $E_2:=E/E_1$.
It is sufficient to show that $\Hom(E,E_2) \cong {\Bbb C}$.
Let $\phi:E \to E_2$ be a non-zero homomorphism.
Then we see that $\ker \phi \cong E_1$ and $\phi$
is surjective.
We assume that $\dim \Hom(E,E_2)>1$.
We set ${\Bbb P}:={\Bbb P}(\Hom(E,E_2)^{\vee})$.
Since $E_1$ is simple, we get an exact sequence
$$
0 \to E_1 \boxtimes {\cal O}_{{\Bbb P}}(n) \to E \to
E_2 \boxtimes {\cal O}_{{\Bbb P}}(1) \to 0,
$$
where $n \in {\Bbb Z}$.
We note that
\begin{align*}
\det p_{{\Bbb P}!}(E_1^{\vee} \otimes E_1 
\boxtimes {\cal O}_{{\Bbb P}}(n)) &\cong {\cal O}_{{\Bbb P}}(2n),\\
\det p_{{\Bbb P}!}(E_1^{\vee} \otimes E_2 
\boxtimes {\cal O}_{{\Bbb P}}(1))&\cong {\cal O}_{{\Bbb P}}(-1).
\end{align*}
Hence we see that 
${\cal O}_{{\Bbb P}} \cong {\cal O}_{{\Bbb P}}(2n)
\otimes {\cal O}_{{\Bbb P}}(-1)$, which is a contradiction.
\end{pf}

By using Lemma \ref{lem:a1} and similar correspondence
$N(v_1,v+v_1,v)$, we see that $\codim_{M(v)^{\mu}}(M(v)^{\mu}_i) \geq 2$
for $i \geq 2$.
Indeed, we get that $\dim M(v)_i^{\mu}=\dim N(v_1,v+v_1,v)_{i-1}-(i-2)
=\dim M(v+v_1)_{i-1}+2 \leq \dim M(v)-2$.
In particular $N(v_1,v,w)$ defines a birational correspondence.
By virtue of Lemma \ref{lem:smooth2}, we also see that
$\codim_{M(v)^{\mu}}(M(v)_i^{\mu}) \geq 3$
for $i \geq 3$.
By using [H2, Cor. 4.7], we obtain the following.

\begin{prop}\label{prop:corres2}
We assume that $M(v)^{\mu} \ne \emptyset$.
Then $M(v)$ is deformation equivalent to $M(w)$, 
in particular $M(v)$ is an irreducible symplectic manifold.
Moreover if $\codim_{M(v)}(M(v) \setminus M(v)^{\mu}) \geq 2$,
then the following diagram is commutative.
\begin{equation*}
\begin{CD}
v^{\perp} @>{R_{v_1}}>> w^{\perp}\\
@V{\theta_v}VV  @VV{\theta_{w}}V\\
H^2(M(v),{\Bbb Z}) @>>{\pi_{w*}\pi_v^*}>
H^2(M(w),{\Bbb Z})
\end{CD}
\end{equation*}
\end{prop}

We shall give an estimate of $\codim_{M(v)}(M(v) \setminus M(v)^{\mu})$,
which depends on $l$ and $\langle v^2 \rangle$.
\begin{lem}\label{lem:mu}
We assume that $\langle v^2 \rangle/2l \geq l$. Then 
$$
\dim(M(v))-\dim(M(v) \setminus M(v)^{\mu})
\geq \frac{\langle v^2 \rangle}{2l}-l+1.
$$
In particular $\codim_{M(v)}(M(v) \setminus M(v)^{\mu}) \geq 2$
for $\langle v^2 \rangle/2l >l$.
\end{lem}

\begin{pf}
Let $E$ be a $\mu$-semi-stable sheaf of
$v(E)=v$ and let
$0 \subset F_1 \subset F_2 \subset \cdots \subset F_t=E$
be a Jprdan-H\"{o}lder filtration of $E$ with respect to $\mu$-stability.
We set $E_i:=F_i/F_{i-1}$.
Then the moduli number of this filtration is bounded by
\begin{align*}
& \sum_{i \leq j} (\dim \Ext^1(E_j,E_i)-\dim \Hom(E_j,E_i))+1\\
\leq & -\sum_{i \leq j}\chi(E_j,E_i)+\sum_{i \leq j}\dim \Ext^2(E_j,E_i)+1\\
\leq & -\chi(E,E)+\sum_{i>j}\chi(E_j,E_i)+
\sum_{i<j}\dim \Ext^2(E_j,E_i)+t+1.
\end{align*}
We set
$v(E):=lr+ldH+a \omega$ and
$v(E_i):=l_i r+l_i dH+a_i \omega$.
Since $\langle v(E_i),v(E_j) \rangle=l_i l_j d^2 H^2-r(l_i a_j+l_j a_i)$,
we see that
\begin{align*}
\sum_{i>j}\chi(E_j,E_i) &= -\sum_{i>j} \langle v(E_j),v(E_i) \rangle\\
&=-\sum_{i>j}\left\{\frac{l_i(l_j d^2 H^2-2r a_j)}{2}
+\frac{l_j(l_i d^2 H^2-2r a_i)}{2}
\right \}\\
&=-\sum_{i>j}\left\{\frac{l_i\langle v(E_j)^2 \rangle}{2l_j}+
\frac{l_j\langle v(E_i)^2 \rangle}{2l_i} \right \}\\
&=-\sum_{i}\frac{(l-l_i)\langle v(E_i)^2 \rangle}{2l_i}.
\end{align*}
We set $\max_i \{l_i\}=(l-k)$.
Let $i_0$ be an integer such that $ \langle v(E_{i_0})^2 \rangle \geq 0$.
Since $\sum_i l_i=l$, we obtain that $t \leq k+1$.  
Since $l-l_i-k \geq 0$ and $\langle v(E_i)^2 \rangle \geq -2$,
 we get that
\begin{align*}
\sum_{i>j} \langle v(E_j),v(E_i) \rangle 
&=k\sum_i \frac{\langle v(E_i)^2 \rangle}{2l_i}+
\sum_i \frac{(l-l_i-k)\langle v(E_i)^2 \rangle}{2l_i}\\
& \geq  k\frac{\langle v(E)^2 \rangle}{2l}-
\sum_{i \ne i_0} (l-l_i-k)\\
& \geq k\frac{\langle v(E)^2 \rangle}{2l}-(l-1-k)k.
\end{align*}
If $r>1$ or $l_i>1$ for some $i$,
then for a general filtration, there are $E_i$ and $E_j$ such that
$\Ext^2(E_j,E_i)=0$.
Therefore we get that 
$\sum_{i<j} \dim \Ext^2(E_j,E_i) \leq (k+1)k/2-1$
for a general filtration.
Then the moduli number of these filtrations is bounded by
\begin{align*}
\langle v^2 \rangle- k\frac{\langle v^2 \rangle}{2l}
+(l-1-k)k+\frac{(k+1)k}{2}-1+t+1 
&\leq 
(\langle v^2 \rangle+2)- \left(k\frac{\langle v^2 \rangle}{2l}
-lk+\frac{k(k-1)}{2}+1 \right)\\
& \leq (\langle v^2 \rangle+2)- \left(\frac{\langle v^2 \rangle}{2l}
-l+1 \right).
\end{align*}
If $l_i=1$ for all $i$ and $r=1$, then $\sum_{i<j}\dim \Ext^2(E_j,E_i) \leq 
k(k+1)/2=(l-1)l/2$, and hence we obtain that
the moduli number is bounded by
$$
(\langle v^2 \rangle+2)- \left(\frac{\langle v^2 \rangle}{2l}
-l+\frac{(l-1)(l-2)}{2} \right).
$$
Therefore we obtain that 
$$
\dim(M(v))-\dim(M(v) \setminus M(v)^{\mu})
\geq \frac{\langle v^2 \rangle}{2l}-l+1
$$
unless $r=1$ and $l=2$.
If $r=1$ and $l=2$, then the moduli number of filtrations
$$
0 \subset T_Z \subset E
$$ of $c_2(I_Z)=c_2(E)$ is 
$(\langle v^2 \rangle+2)-(\langle v^2 \rangle/2l-l)$.
In this case, we can easily show that 
$E^{\vee \vee} \cong {\cal O}_X^{\oplus 2}$,
and hence the moduli number of these $E$ is bounded by
$3 c_2(E)-3 \leq (\langle v^2 \rangle+2)-(\langle v^2 \rangle/2l-l+1)$.
Therefore we obtain our lemma. 
\end{pf}

\begin{thm}\label{thm:2}
Let $v=lr+l \xi+a \omega, \;\xi \in H^2(X,\Bbb Z)$
 be a primitive Mukai vector such that $r>0$ and $r+\xi$ is primitive.
If $\langle v^2 \rangle/2 > 
\max\{rl(rl-1),l^2\}$, then 
$M(v)$ is deformation equivalent to 
$\Hilb_X^{\langle v^2 \rangle/2+1}$ and
$$
\theta_v : v^{\perp} \longrightarrow H^2(M(v),{\Bbb Z})
$$
is an isometry which preserves hodge structures.
\end{thm}

\begin{pf}
By Proposition \ref{prop:deform} and Lemma \ref{lem:mu},
it is sufficient to find a Mukai vector 
$v'=l(r+\xi')+a' \omega$ which satisfies
(1) $r+\xi'$ is primitive, (2) $\langle {v'}^2 \rangle=
\langle {v}^2 \rangle$, (3) $a' \equiv a \mod l$, and 
(4) $M(v')$ satisfies our claims.
By Proposition \ref{prop:corres2},
it is sufficient to find a Mukai vector
$v'$ of the form \eqref{eq:4-1}.
Hence $a \equiv -p \mod l$.
We have to choose integers $r_1,d_1$ and $d$.
We note that
$pr_1 \equiv -1 \mod l$, and 
$r_1$ and $r$ are relatively prime.
We can easily choose such an integer $r_1$ of 
$0<r_1<lr$.
Then we can choose $d$ and $d_1$ of $d r_1-d_1r=1$.
Since $\langle v^2 \rangle/2l > 
r(rl-1)$, we see that 
$k(s)=(r_1 \langle v^2 \rangle/2l+r)r_1/l-r^2>0$.
Hence we get a Mukai vector
$v'$ of the form \eqref{eq:4-1} 
 \end{pf}
Tohru Nakashima taught the author the following remark.
\begin{rem}\label{rem:exist}
If $r=1$ and $\xi=0$, then $rl(rl-1)<l^2$.
Hence the assertions of Theorem \ref{thm:2} hold for
$\langle v^2 \rangle/2l > l$.
It is known that 
$M(v)^{\mu}\ne \emptyset$ if and only if
$\langle v^2 \rangle/2l \geq l$,
and $\langle v^2 \rangle/2l \ne l$ if $v$ is primitive.
\end{rem}

\vspace{1pc}

{\it Acknowledgement.}
The author would like to thank Masaki Maruyama for
giving him a preprint of Ellingsrud and Str\o mme ([E-S]),
which inspired him very much,
Tohru Nakashima for teaching him Remark \ref{rem:exist},
and Hiraku Nakajima for valuable discussions.


\begin{thebibliography}{[G-H]}




\bibitem[B]{B:1}
Beauville, A.,
{\it Vari\'{e}t\'{e}s K\"{a}hleriennes dont la premi\`{e}re classe Chern
est nulle,}
J. Diff. Geom. {\bf 18} (1983), 755--782
\bibitem[D1]{D:1}
Drezet, J.-M.,
{\it Fibr\'{e}s exceptionnels et vari\'{e}t\'{e}s de
modules de faisceaux semi-stables sur ${\Bbb P}_2({\Bbb C})$,}
J. reine angew. Math. {\bf 380} (1987), 14--58
\bibitem[D2]{D:2}
Drezet, J.-M.,
{\it Cohomologie des vari\'{e}t\'{e}s de modules
de hauteur nulle,}
Math. Ann. {\bf 281} (1988), 43--85
\bibitem[D-L]{D-L:1}
Drezet, J.-M., Le-Potier, J.,
{\it Fibr\'{e}s stables et fibr\'{e}s exceptionnels sur ${\Bbb P}^2$,}
Ann.\ scient.\ \'{E}c. Norm. Sup., $4^e$ s\'{e}rie, t. {\bf 18} (1985),  193--244   
\bibitem[E-S]{E-S:1} 
Ellingsrud, G., Str\o mme, S. A.,
{\it Towards the Chow ring of the Hilbert scheme on ${\Bbb P}^2$,}
J. reine angew. Math. {\bf 441}
(1993), 33--44
\bibitem[G-H]{G-H:1}
G\"{o}ttsche, L., Huybrechts, D.,
{\it Hodge numbers of moduli spaces of stable bundles on K3 surfaces,}
Internat. J. Math. {\bf 7} (1996), 359--372
\bibitem[H1]{H:1}
Huybrechts, D.,
{\it Birational symplectic manifolds and their deformations,}
J. Diff. Geom. {\bf 45} to appear
\bibitem[H2]{H:2}
Huybrechts, D.,
{\it Compact Hyperk\"{a}hler Manifolds: Basic Results,}
alg-geom/9705025
\bibitem[Ma1]{M:1}
 Maruyama, M.,
 {\it Moduli of stable sheaves II,}
 J. Math. Kyoto Univ. {\bf 18} (1978),  557--614
\bibitem[Ma2]{M:3}
 Maruyama, M.,
{\it Moduli of algebraic vector bundles,}
in preparation
\bibitem[Mu1]{Mu:2}
Mukai, S.,
{\it Symplectic structure of the moduli space of sheaves on an 
abelian or K3 surface,}
Invent. math. {\bf 77}
(1984), 101--116
\bibitem[Mu2]{Mu:3}
Mukai, S.,
{\it On the moduli space of bundles on K3 surfaces I,}
Vector bundles on Algebraic Varieties, Oxford, 1987, 341--413 
\bibitem[Mu3]{Mu:5}
Mukai, S.,
{\it Moduli of vector bundles on K3 surfaces, and symplectic manifolds,}
Sugaku Expositions, {\bf 1} (1988), 139--174 
\bibitem[O]{O:1}
O'Grady, K.,
{\it The weight-two Hodge structure of moduli spaces of sheaves
on a K3 surface,}
preprint (1995)
\bibitem[S]{S:1}
Simpson, C.,
{\it Moduli of representations of the fundamental group
of a smooth projective variety I,}
Publ. Math. I.H.E.S. {\bf 79} (1994), 47--129
\bibitem[Y]{Y:1}
Yoshioka, K.,
{\it Some notes on the moduli of stable sheaves on elliptic surfaces,}
alg-geom/9705007   
   
\end{thebibliography}
\end{document}